\documentclass[12pt,reqno,tbtags]{amsart}
\usepackage{amssymb}
\numberwithin{equation}{section}

\usepackage[usenames]{color}

\newcommand{\halmos}{\rule{1ex}{1.4ex}}
\newcommand{\proofbox}{\hspace*{\fill}\mbox{$\halmos$}}
\newenvironment{proofof}[1]{\noindent {\it
Proof of #1}.}{\proofbox\par\smallskip\par}
\newenvironment{prf}{\noindent {\it
Proof}.}{\proofbox\par\smallskip\par}

\newtheorem{theorem}{Theorem}[section]
\newtheorem{lemma}[theorem]{Lemma}

\theoremstyle{definition}

\theoremstyle{remark}

\newcounter{thmenumerate}

\newcounter{xenumerate}

\newcommand\REM[1]{\texttt{[#1]}\marginal{XXX}}

\begingroup
  \divide\count255 by 60
  \multiply\count255 by -60
  \advance\count255 by \time
  \ifnum \count255 < 10 \xdef\klockan{\the\count1.0\the\count255}
  \else\xdef\klockan{\the\count1.\the\count255}\fi
\endgroup

\def\rompar(#1){\textup(#1\textup)}    

\def\xexp(#1){e^{#1}}

\newcommand\op{o_{\mathrm p}}

\newcommand\bbN{\mathbb N}

\newcommand\E{\operatorname{\mathbb E{}}}

\newcommand\Bi{\operatorname{Bi}}

\newcommand\eps{\varepsilon}

\renewcommand\phi{\varphi}

\def\[#1]{[\![#1]\!]}

\newcommand\whp{\textbf{whp}}

\newcommand{\ind}{\mathbb{I}}

\newcommand\marginal[1]{\marginpar{\raggedright\parindent=0pt\tiny #1}}

\newcommand\REMRM[1]{}

\newcommand{\eq}{\begin{equation}}
\newcommand{\en}{\end{equation}}
\newcommand{\lbeq}[1]{\label{#1}}
\newcommand{\refeq}[1]{(\ref{#1})}

\newcommand{\eqn}[1]{\eq #1 \en}
\newcommand{\cn}{\Omega}
\newcommand{\vep}{\varepsilon}
\newcommand{\prob}{{\mathbb P}}
\newcommand{\Qprob}{{\mathbb Q}}

\newcommand{\sss}{\scriptscriptstyle}

\newcommand{\eqan}[1]
{\begin{align}
#1
\end{align}
}

\renewcommand{\Pr}{{\mathbb P}_p}

\newcommand{\een}{\end{enumerate}}

\newcommand{\ben}{\begin{enumerate}}
\newcommand{\beq}{\begin{equation}}
\newcommand{\eeq}{\end{equation}}
\newcommand{\bec}{\begin{center}}
\newcommand{\ece}{\end{center}}

\begin{document}
\title{The second largest component in the supercritical 2D Hamming graph}

\date{3 January 2009 (\today)}

\author{Remco van der Hofstad}
\address{Department of Mathematics and
        Computer Science, Eindhoven University of Technology,
        5600 MB Eindhoven, The Netherlands.}
\email{rhofstad@win.tue.nl}
\urladdr{http://www.win.tue.nl/$\sim$rhofstad}

\author{Malwina J. Luczak}
\address{Department of Mathematics, London School of Economics,
  Houghton Street, London WC2A 2AE, United Kingdom}
\email{m.j.luczak@lse.ac.uk}
\urladdr{http://www.lse.ac.uk/people/m.j.luczak@lse.ac.uk/}

\author{Joel Spencer}
\address{Computer Science and Mathematics Departments, New York
University, Mercer Street,
         New York City, United States of America.}
\email{spencer@courant.nyu.edu}
\urladdr{http://cs.nyu.edu/cs/faculty/spencer/}

\keywords{random graphs, percolation, phase transition, scaling window}
\subjclass[2000]{05C80}

\begin{abstract}
The $2$-dimensional Hamming graph $H(2,n)$ consists of the $n^2$
vertices $(i,j)$, $1\leq i,j\leq n$, two vertices being adjacent
when they share a common coordinate.   We examine random
subgraphs of $H(2,n)$ in percolation with edge probability $p$, in
such a way that the average
degree satisfies $2(n-1)p=1+\varepsilon$.  Previous work~\cite{hl} has shown that
in the barely supercritical region $n^{-2/3}\ln^{1/3}n\ll \varepsilon
\ll 1$, the largest component satisfies a law of large numbers
with mean $2\varepsilon n$. Here we
show that the second largest component has, with high probability, size bounded by
$2^8\varepsilon^{-2}\log(n^2\varepsilon^3)$, so that the dominant component
has emerged. This result also suggests that a {\it discrete duality
principle} holds, where, after removing the largest connected component
in the supercritical regime, the remaining random subgraphs behave as in
the subcritical regime.
\end{abstract}

\maketitle

\section{Introduction and main result}\label{S:intro}

In their seminal work~\cite{er}, Paul Erd\H{o}s and Alfred R\'enyi
noted with surprise the development of a {\em giant component} in
the random graph $G(n,p)$ where each of the $n(n-1)/2$ possible edges
of the complete graph of size $n$ is present with probability $p$ independently
of all the other edges.  When the average degree $(n-1)p$ satisfies
$(n-1)p=1+\vep$, and $\vep$ is positive and fixed (independent of
$n$), then the largest component will contain a positive proportion
of the vertices while the size of the second largest component is
only logarithmic in $n$. Today we see this as a phase transition
phenomenon exhibiting what mathematical physicists call `mean-field'
behaviour.

For many years, there has been great interest in the {\em barely
supercritical} phase of $G(n,p)$, that is, the range of $p$ values
given by $p=(1+\vep)/n$, where $\vep=\vep(n)$ satisfies
$n^{-1/3}\ll \vep = o(1)$.  For convenience, we can also write
$\vep=\lambda n^{-1/3}$, where $\lambda=\lambda(n)\rightarrow +\infty$, but does
so more slowly than $n^{1/3}$. In this phase the {\em dominant
component} has already appeared. We actually know quite precisely
that the largest component, $C_{\sss (1)}$, satisfies $|C_{\sss (1)}|=
2\vep n(1+\op(1))=2\lambda n^{2/3}(1+\op(1))$ with probability tending to 1 as $n \to
\infty$, where $\op(1)$ denotes a quantity that converges to
zero in probability; and that the second largest component, $C_{\sss (2)}$, satisfies
$|C_{\sss (2)}| = \Theta(\vep^{-2}\log(n\vep^3))=
\Theta(n^{2/3}\lambda^{-2}\ln\lambda)$ with probability tending to 1
as $n \to \infty$. Thus, in particular, $|C_{\sss (2)}|\ll n^{2/3} \ll
|C_{\sss (1)}|$. Further, as $\lambda$ increases, the largest component
increases in size while the second largest diminishes in size.
(Actually, the second largest component is being frequently
`gobbled up' by the dominant component, generally leaving a
smaller component as the new second component.)
See \cite{Boll84b, Lucz90a} for the proofs of these results,
and \cite{Boll01, JLR} for introductions to the field.
We feel, speaking quite generally, that an intensive
study of the second largest component is vital to enhancing our
understanding of percolation phenomena.

We believe that the second largest component should grow until the random
structure reaches a {\em critical window}.  In that critical window,
which, for $G(n,p)$, means $p=(1+\lambda n^{-1/3})/n$ with
$\lambda$ fixed, the first and second
largest components exhibit complex chaotic behaviour. On the other hand, in the
barely supercritical phase, just after the critical window, the
dominant, or `giant', component will have asserted itself.


In the present paper, our object of study is the $2$-dimensional
Hamming graph $H(2,n)$.  The $n^2$ vertices of this graph can be
represented as ordered pairs $(i,j)$, $1\leq i,j \leq n$.  Vertices
$(i,j)$ and $(i',j')$ are adjacent if and only if either $i=i'$ or $j=j'$.
Pictorially, $H(2,n)$ consists of an $n\times n$ lattice with each
horizontal and vertical line being a complete graph. We write $\cn
=2(n-1)$ for the vertex degree in $H(2,n)$. We examine random
subgraphs of $H(2,n)$ in independent percolation with edge probability
$p$; that is, each edge is kept with probability $p$ and removed
with probability $1-p$, independently of all other edges. We set
$p_c=1/\cn$, which will act as our critical probability;
a justification for this definition of critical probability lies in
the recent results in~\cite{bchss1,bchss2,hl}. We parametrise $p=
(1+\vep)/\cn$ so that the average vertex degree is $1+\vep$.
Throughout the rest of the paper $C_{\sss (1)}, C_{\sss (2)}$ will refer to the
largest and second largest components, respectively, of the Hamming
graph $H(2,n)$. Also, we shall use the phrase `with high probability'
(\whp{}) to mean `with probability $1-o(1)$ as $n \to \infty$'.

This work continues the exploration of van der Hofstad and
Luczak~\cite{hl}. It was shown therein that, when $n^{-2/3}(\log
n)^{1/3} \ll \vep \ll 1$, the largest component has size $2\vep n^2
(1+\op(1))$.   The general sense of a mean-field percolation event in
percolation on a graph with $V$ vertices is that there is a critical
probability $p_c$, and that the barely supercritical phase occurs
when $p=p_c(1+\vep)$ and $V^{-1/3} \ll \vep \ll 1$. This is the case
in the Erd\H{o}s-R\'enyi phase transition with $V=n$. For the
$H(2,n)$ phase transition, $V=n^2$, and so the above results, up to a
logarithmic term, fit the mean-field paradigm. Here we study the
second largest component in percolation on $H(2,n)$ in the barely
supercritical region.  In this aspect we are also able to, again up to
a logarithmic term, fit the mean-field paradigm. In the mean-field
picture of random graphs, the structure remaining when the dominant
component is removed for $p=p_c(1+\vep)$ where $\vep=\vep(n)
\rightarrow 0$  is like  the largest connected component in
the subcritical regime with $p=p_c(1-\vep)$. It is well known
that in this regime in $G(n,p)$, the second largest component
is of order $\vep^{-2} \log (n\varepsilon^3)$. The upper bound
is the content of our main result. In our results
for $H(2,n)$, we shall always work
at $p=p_c+\vep/\cn=p_c+\vep/(2(n-1))$. In \cite[(1.10) and (1.11)]{hl},
it is shown that there is little difference in working
with $p=p_c+\vep/\cn$ or $p=(1+\vep)/\cn$, and we refer the
reader there for more details.

\begin{theorem}
[The second component in the supercritical phase for $H(2,n)$]
\label{thm-main} Consider the 2-dimensional Hamming graph 
$H(2,n)$. Let $p=p_c +\frac{\vep}{\cn}$ and let
$n^{-2/3}(\log{n})^{1/3}\ll \vep \ll 1$. Then, \whp{},
    \eqn{
    |C_{\sss (2)}| \le 2^8 \vep^{-2} \log (n^2\varepsilon^3).
    }
\end{theorem}

In particular, this result implies that the ratio of the sizes of
the second and first largest components tends to zero in this regime, a salient
feature of the barely supercritical phase.  We feel that this
feature should hold even without the logarithmic separation from
criticality. That is, parametrise $2(n-1)p=1+\varepsilon$ and assume
only $n^{-2/3}\ll \varepsilon \ll 1$. We conjecture,
following~\cite{hl}, that the largest component will have size
$2\varepsilon n^2(1+\op(1))$.  We further conjecture that the second largest
component will have size $\ll n^{2/3}$, which in particular is
asymptotically smaller than the largest component. Let us note at this point that the
logarithmic gap from the critical window (defined as
in~\cite{bchss1,bchss2,bchss3}) has recently been removed by Asaf
Nachmias~\cite{n07}; however, he does not establish a law of large
numbers for the giant component, and he does not consider the second
largest component. Further, we conjecture that, when $n^{-2/3}\ll \varepsilon \ll 1$,
$|C_{\sss (2)}|=\Theta(\vep^{-2} \log (n^2\varepsilon^3))$, i.e.,
the bound in Theorem \ref{thm-main} is sharp. Thus, in particular, we conjecture
that the barely supercritical regime for $H(n,2)$ has similar behaviour
to that of $G(n,p)$.

\section{Preliminaries}

In this section, we establish a lemma for a class of branching processes that
will play a key role in our proofs.

We start with an inequality concerning deviations of binomial
random variables below their mean.
%
%
If $X\sim \Bi(k,p)$, then
(see for instance~\cite{Jans02})
    \eqn{
    \lbeq{Binbd}
    \prob(X\leq kp-t) \leq e^{-\frac{t^2}{2(kp+\frac{t}{3})}}.
}

We consider Galton-Watson processes where each individual's
offspring is a random variable $Z$ such that $\E[Z^2] < \infty$. We
always assume that our process begins with one individual. Sometimes
we shall take $Z$ to have a binomial distribution $\Bi (N,p)$, with
$p$ the Hamming graph edge probability, and $N$ a suitable positive
integer. We will write $\prob_{N,p}$ for the probability measure
corresponding to this process.
We will also need Galton-Watson processes that are
`inhomogeneous', in that the offspring size may vary depending on
the parent's `location' in the Galton-Watson tree.
\medskip

A Galton-Watson process can be thought of as a 2-dimensional Markov
chain $(Q_t,G_t)$, where $Q_t$ is the total progeny born until time
$t$, and $G_t$ is the total number of `active' population members,
that is those that are yet to have offspring. To be precise, we
think of a Galton-Watson process as an evolving tree that is
explored one node at a time; then $Q_t$ is the total number of nodes
in that tree at time $t$, and $G_t$ is the total number of
unexplored nodes at time $t$. At each time $t$, if $G_t
> 0$, then we choose one active member of the population and decide
the number of its offspring. In a homogeneous Galton-Watson process,
all population members have the same offspring distribution, in our
case binomial distribution $\Bi (N,p)$, where $N \in \bbN$ and $p
\in [0,1]$ is the Hamming graph edge probability. Then, on the event
$\{ \omega: G_t (\omega) \ge 1\}$,
    \eqn{
    Q_{t+1}= Q_t + Z_t,\qquad
    G_{t+1}= G_t + Z_t-1,
    }
where the $Z_t$ are i.i.d.\ $\Bi (N,p)$. We always assume $Q_0=G_0=1$; we
let $\phi_0 = \varnothing$ and
    \eqn{
    \phi_t = \sigma (Z_s: s \le t), \quad \quad \quad t=1,2, \ldots,
    }
the $\sigma$-field generated by $Z_s$ ($s=1,2, \ldots,t$).

Note that, by the above,
    \eqan{
    \ind_{G_{t-1} \ge 1} Q_t &= \ind_{G_{t-1} \ge 1}(1 + Z_0 + \ldots + Z_{t-1}),\\
    \ind_{G_{t-1} \ge 1} G_t &= \ind_{G_{t-1} \ge 1}(1+ Z_0 + \ldots + Z_{t-1}-t),
    }
where $\ind_{\mathcal{A}}$ denotes the indicator of the event $\mathcal{A}$.
Letting $T_0 = \inf \{t: G_t =0\}$, we further see that
    \eqan{
    Z_0 + \ldots + Z_{T_0-1}=T_0-1.
    }

For convenience, we shall instead assume that our Galton-Watson
processes have a continuation in that the random variables $Z_t$
continue to be generated even after $G_t$ has hit 0. With this
assumption, we may simply write, for each $t$,
    \eqn{
    Q_t = 1 + Z_0 + \ldots + Z_{t-1},\qquad  G_t = 1+ Z_0 + \ldots + Z_{t-1}-t.
    }

We shall now generalise this setting to a branching
process where $Z_t$ are not i.i.d., but where each $Z_t$ is Binomial
$\Bi (N_t,p)$. Here, $p$ is the edge probability of the Hamming
graph $H(2,n)$, and each $N_t$ is a random variable independent of
$Z_t$ and such that, for each $t$,
and conditionally on $\phi_{t-1}$,
    \begin{equation}
    \label{narrow-band} N_t p \in [1+\vep/2,1+\vep]  \quad \quad
    \mbox{with probability } 1.
    \end{equation}
Any such branching process will be called generically a {\it
narrow-banded (discrete) birth process}, and in the lemma below we
use $\prob$ to denote the corresponding probability measure -- with
a slight abuse of notation, since in fact there is a whole class of
processes satisfying property~(\ref{narrow-band}).

\begin{lemma}[Large finite progeny for narrow-banded processes]
\label{lem.key}
Let $Q$ be the
total population size of a narrow-banded birth process as
defined above.
Then, for $\alpha > 0$ sufficiently large, there exists a constant $C$ such that
    \begin{eqnarray*}
    \prob (\alpha \vep^{-2} \le Q <\infty) < C \vep e^{-\alpha/2^8}.
    \end{eqnarray*}
\end{lemma}

\begin{prf}
As described above, at each time $t$ we generate $Z_t$, a binomial
$\Bi (N_t,p)$, where
    \eqn{
    p=\frac{1+ \vep}{2(n-1)}\quad \text{and}\quad \prob (N_t p\in [1+ \vep/2,1 + \vep]  \mid \phi_{t-1})=1.
    }

As earlier in this section, let $G_0=1$, and let
$G_t=1+Z_0+\ldots+Z_{t-1}-t$ for $t=1,2, \ldots$. In other words, we
assume our narrow-banded Galton-Watson process to have a
continuation; clearly, this does not in any way affect the correctness of the proof
of the lemma.

Let
    \eqan{
    {\mathcal E}&=\Big \{G_{\lfloor \alpha \vep^{-2}/2 \rfloor } < \alpha
    \vep^{-1}/16; \quad G_t > 0 \mbox{ for all } t\leq \alpha
    \vep^{-2}/2\Big \},\\
    {\tilde {\mathcal E}}&=\{G_{\lfloor \alpha \vep^{-2}/2 \rfloor } \geq \alpha
    \vep^{-1}/16; \quad G_t > 0 \mbox{ for all } t\leq \alpha
    \vep^{-2}/2\}.
}
Let $T_0$ be the first time $t$ such that $G_t=0$. Then
$T_0 \le \alpha \vep^{-2}/2$, and $G_{T_0-1}=1$, so it follows that
    \begin{eqnarray*}
    Q = Q_{T_0-1}  = G_{T_0-1} + T_0-1  \le  \alpha \vep^{-2}/2 < \alpha
    \vep^{-2}.
    \end{eqnarray*}
Therefore, we can upper bound
    \eqan{
    \prob (\alpha \vep^{-2} \le Q <\infty)&=\prob (\{\alpha \vep^{-2} \le Q <\infty\}\cap  {\mathcal E}^c)
    +\prob (\{\alpha \vep^{-2} \le Q <\infty\}\cap  {\mathcal E})\nonumber\\
    &\leq \prob (\{\alpha \vep^{-2} \le Q <\infty\}\cap {\tilde
    {\mathcal E}})+\prob({\mathcal E}) \nonumber\\
    & =\prob (\alpha \vep^{-2} \le Q <\infty\mid {\tilde {\mathcal E}})
    \prob ({\tilde {\mathcal E}}) +\prob({\mathcal E}),
    }
since, if $Q\geq \alpha \vep^{-2}$, then $G_t>0$ for all $t\leq \alpha \vep^{-2}/2$.
We shall start by bounding $\prob (\alpha \vep^{-2} \le Q
<\infty\mid {\tilde {\mathcal E}})$, and later bound
$\prob({\mathcal E})$ and $\prob ({\tilde {\mathcal E}})$.

%
%

On the event $\tilde{\mathcal E}$, for $t \ge \alpha
\vep^{-2}/2$, we couple our narrow-banded birth process with the
standard Galton-Watson process where the offspring distribution is
binomial
    \eqn{
    \lbeq{Bin-law-min}
    \Bi \Big (\Big \lfloor \frac{2(n-1)(1+ \vep/2)}{1+ \vep}\Big \rfloor, \frac{1+
    \vep}{2(n-1)} \Big ),
    }
with mean less than or equal to $1+\vep/2$. Let
$Q^{-}_t$ and $G^{-}_t$ respectively be the total progeny and total
number of active members at time $t$ for this process; assume also
that $Q^{-}_{\lfloor \alpha \vep^{-2}/2\rfloor} =\lceil
\alpha \vep^{-1}/16 \rceil$. Let $t_0 = \lfloor \alpha
\vep^{-2}/2 \rfloor $, and for $t \ge t_0+1$ let us write
    \eqn{
    \lbeq{QGt-def-min}
    Q^-_t= \lceil \alpha \vep^{-1}/16 \rceil +Z^-_{t_0} +
    \ldots +
    Z^-_{t},\qquad
    G^-_t = \lceil \alpha \vep^{-1}/16 \rceil + Z^-_{t_0} +
    \ldots + Z^-_{t}-(t-t_0),}
where the $Z^-_t$ are i.i.d.\ $\Bi\big(
\lfloor (2(n-1) (1+ \vep/2)/(1+ \vep) \rfloor, (1+ \vep)/2(n-1)\big)$.
Let $\prob^-$ denote the probability measure corresponding to this
process.

The coupling is between the corresponding tree exploration
processes, step-by-step, as is standard (and as used, for instance, in~\cite[Section 4]{hl}), so that
$Z_t \ge Z^-_t$, $Q_t \ge Q^{-}_t$ and $G_t \ge G^{-}_t$ for all $t
\ge \lfloor \alpha \vep^{-2}/2 \rfloor$. This implies that
    \eqn{
    \lbeq{prob-cond-Etilde}
    \prob (\sup_t Q_t < \infty\mid \tilde{\mathcal E} ) \le
    \prob^- (\sup Q^{-}_t < \infty).
    }
Let $\Qprob_k$ denote the law of a branching process with offspring distribution
in \refeq{Bin-law-min}, and starting from $k$ individuals, and let
$\tilde Q^-$ be its total progeny. Then,
    \eqn{
    \prob^- (\sup Q^{-}_t < \infty)= \Qprob_{\sss \lceil \alpha \vep^{-1}/16 \rceil}(\tilde Q^-<\infty)
    =\Qprob_{1}(\tilde Q^-<\infty)^{\lceil \alpha \vep^{-1}/16 \rceil},
    }
by the independence of the evolution of the initial individuals.

Further, \cite[(3.28)]{hl} shows that
    \eqn{
    \Qprob_{1}(\tilde Q^-<\infty)=1-\vep +O(n^{-1}+\vep^2),
    }
so that
    \eqn{
    \lbeq{prob-min-bd}
    \prob^- (\sup
    Q^{-}_t < \infty) = (1-\vep + O(n^{-1}+ \eps^2))^{\lceil \alpha
    \vep^{-1}/16 \rceil } \le e^{-\alpha /16}(1+O(\alpha (\vep + (\vep n)^{-1}))),
    }
for all $\vep \in (0,1)$ and $\alpha >0$.
We conclude from the above thet
    \begin{eqnarray*}
\prob (\sup_t Q_t < \infty\mid \tilde{\mathcal E} )
    \le  e^{-\alpha /16}(1+o(1)).
    \end{eqnarray*}


Next we shall show that $\prob({\mathcal E})$ and $\prob ({\tilde
{\mathcal E}})$ are quite small.
For this, we also need an upper bounding Galton-Watson process
$(Q^+_t,G^+_t)$, where
the offspring distribution is binomial
$$\Bi \Big (2(n-1), \frac{1+ \vep}{2(n-1)} \Big ),$$
with mean $1+\vep$.  Let $Q^{+}_t$ and $G^{+}_t$ respectively be the
total progeny and total number of active members at time $t$ for this process;
assume also that the initial population size is $Q^{+}_0 =1$. Let us write
    \eqn{
    Q^+_t= 1 +Z^+_0 + \ldots +
    Z^+_{t-1},\qquad
    G^+_t = 1 + Z^+_0 + \ldots + Z^+_{t-1}-t,
    }
where the $Z^+_t$ are
i.i.d.\ $\Bi (2(n-1), (1+ \vep)/2(n-1))$. We use $\prob^+$ to denote the
corresponding probability measure.

Similarly, we use a lower bounding Galton-Watson process
$(Q^-_t,G^-_t)$ (but this time starting from time 0, rather than
from time $t_0 = \lfloor \alpha \vep^{-2}/2 \rfloor $), where the
offspring distribution is binomial $\Bi \big(\lfloor 2(n-1)(1+
\vep/2)/(1+ \vep) \rfloor, (1+ \vep)/2(n-1)\big)$, with mean at most
$1+\vep/2$. Let $Q^{-}_t$ and $G^{-}_t$ be the total progeny and
total number of active members at time $t$ for this process; assume
also that the initial population size is $Q^{-}_0 =1$. Let us write
    \eqn{
    Q^-_t= 1 +Z^-_0 + \ldots +
    Z^-_{t-1},\qquad
    G^-_t= 1 + Z^-_0 + \ldots + Z^-_{t-1}-t,}
where the $Z^-_t$ are
i.i.d.\ $\Bi \big(\lfloor 2(n-1)(1+ \vep/2)/(1+ \vep) \rfloor, (1+
\vep)/2(n-1)\big)$. Once again, we use $\prob^-$ to denote the
corresponding probability measure.

Now we couple $Z_t$ with $Z^-_t$ and $Z^+_t$, so that, for all $t = 0,1,
\ldots$,
    \eqn{
    Z^-_t \le Z_t \le Z^+_t.
    }
A suitable coupling can be achieved thanks to standard results about
stochastic domination between binomial random variables with different
parameters. Explicitly, for $n$ large enough, we may generate
independent binomial random variables $Z^-_t$ such that
    $$
    Z^-_t \sim \Bi
    \Big(\lfloor 2(n-1)(1+ \vep/2)/(1+ \vep) \rfloor, (1+ \vep)/2(n-1)\Big),
    $$
and independent binomial
random variables $W_t$ such that
    $$
    W_t \sim \Bi \Big(2(n-1) -\lfloor 2(n-1)(1+ \vep/2)/(1+ \vep)\rfloor,
    (1+ \vep)/2(n-1)\Big).
    $$
We can then set $Z^+_t = Z^-_t + W_t$ for all $t$. We let
$\prob^{+,-}$ denote the coupling measure.

Let ${\mathcal A}$ be the event that $G^+_t > 0$ for all $t \le
t_0$.
Let ${\mathcal B}$ be the event that $G^-_{t_0} <
\alpha \eps^{-1}/16$.
Note that, under the coupling,
    \eqn{
    {\mathcal E} \subseteq {\mathcal A} \cap {\mathcal B}.
    }

However, now it is easily seen
(using \cite[Proposition 3.2]{hl}) that there exists $c_0= c_0
(\alpha)$ such that $c_0 \to 0$ as $\alpha \to \infty$ and
    \eqn{
    \Pr^{+,-}
    ({\mathcal A}) =\Pr(Q\geq t_0) =2 \vep+O\big(\vep^2 +1/\sqrt{t_0}\big)
    \le (2+c_0(\alpha)) \vep (1+o(1)).
    }
This follows since $t_0=\lfloor \alpha\vep^{-2}/2 \rfloor$. Also, using~(\ref{Binbd}),
    \eqan{ \Pr^{+,-} ({\mathcal B}) & \le \Pr
    \Big (\Bi \big(\lfloor 2(n-1)(1+ \vep/2)/(1+\vep) \rfloor \lfloor \alpha
    \vep^{-2}/2 \rfloor, (1+ \vep)/2(n-1)\big) < \alpha \vep^{-2}/2 +
    \alpha \vep^{-1}/16 \Big )\nonumber\\
    & \le \exp \big (-\alpha/2^8\big),
    }
for all $\vep$ satisfying $\vep \gg n^{-2/3} (\log n)^{1/3}$, and all $\alpha$ and
$n$ sufficiently large, since $\alpha \vep^{-1} \gg 1$ and $\alpha
\vep^{-1} \gg \alpha \vep^{-2}/n$.

Now, the event ${\mathcal A}$ is increasing, and the event ${\mathcal B}$ is
decreasing, and both are events on the same probability space,
corresponding to a family of independent random variables. It then
follows from the FKG inequality that they are negatively correlated.
Hence,
    \eqn{ \Pr ({\mathcal E} ) \leq \Pr^{+,-} ({\mathcal A} \cap
    {\mathcal B}) \leq \Pr^{+,-}  ({\mathcal A}) \Pr^{+,-}( {\mathcal
    B}) \le (2+c_0(\alpha)) \vep e^{-\alpha /2^8}(1+o(1)).
    }
Also,
    $$\prob ({\tilde {\mathcal E}}) \le \Pr^{+,-} ({\mathcal A}),$$
and hence
    \eqn{
    \Pr (\alpha \vep^{-2} \le Q < \infty) \le \prob (\{\alpha \vep^{-2} \le Q <\infty\}\cap {\tilde {\mathcal E}})+\prob({\mathcal E})
    \le 2(2+c_0 (\alpha)) \vep e^{-\alpha /2^8} (1+o(1)).
    }
\end{prf}

\section{Proof of main result}

Recall that $\Omega = 2(n-1)$. Let $Q({\bf
v})$ denote the component of vertex ${\bf v}$. Our first lemma is
\cite[Proposition 2.1]{hl}.
\begin{lemma}[Cluster tail equals the survival probability]
\label{lem.large-clusters}
Let $\vep$ satisfy $\vep^3 n^2 \gg \log n$. Let $p=p_c
+\vep/\cn$. Let $N \gg \vep^{-2}$. Then, for any vertex ${\bf
v}_0=(i_0,j_0)$,
    \eqn{
    \Pr \Big ( |Q({\bf v}_0)| \ge N  \Big ) = 2 \vep (1+o(1)).
    }
\end{lemma}

Our next lemma upper bounds the variance of $Z_{\sss \ge N}$, the number of vertices
in components at least $N$. This result is a
special case of \cite[Corollary 2.3]{hl}.

\begin{lemma}[Concentration of vertices in large clusters]
\label{lem.variance} Let $\vep$ satisfy $\vep^3 n^2 \gg \log n$.
Let $p=p_c+\vep/\cn$. Let $N \gg \vep^{-2}$. Then, for every
$\delta > 0$,
    \eqn{
    \Pr (|Z_{\sss \ge N} - \E_p [Z_{\sss \ge N}]| \ge \delta \vep n^2)=o(1).
    }
\end{lemma}

We now show that, \whp{}, there are no components of
`medium' size; that is, if $\alpha = \alpha (n) \to \infty$ as $n
\to \infty$, then any component of size at least
 $\alpha\vep^{-2}$, \whp{} will in fact be of size at least
$\varepsilon n^2/5$. This is the content of our next lemma:


\begin{lemma}[No middle ground]
\label{lem.no-middle-ground}
Let $\vep$ satisfy $\vep^3 n^2 \gg \log n$. Let $p=p_c
+\vep/\cn$. Let $\alpha = \alpha (n) \to \infty$ as $n \to
\infty$. Then there exists a constant $C$ such that, for $n$ large
enough,
    \eqn{
    \lbeq{no-middle-ground}
    \Pr \Big ( \alpha \vep^{-2}  \le |Q({\bf v_0})| < \frac{ \varepsilon
    n^2}{5}\Big ) \le C \big(\vep e^{-\alpha/2^8} +n^{-6}\big).
    }
Hence, the probability that there is some vertex ${\bf v_0}$ such
that its component $Q({\bf v_0})$ satisfies
    \eqn{
    \lbeq{no-comp}
    2^8 \vep^{-2} \log
    (n^2 \vep^3) \le |Q({\bf v_0})| < \vep n^2/5
    }
is $o(1)$ as $n \to \infty$.
\end{lemma}

Before giving the proof of Lemma~\ref{lem.no-middle-ground}, let us
state two more results from~\cite{hl}, which compare the size of the
cluster of a vertex to the total progeny of suitable Galton-Watson
processes.

The first of these is essentially \cite[Lemma 4.1]{hl}, proved
by standard methods, and gives an upper bound. Let $\prob_{\cn,p}$
be the probability measure corresponding to a standard Galton-Watson
process where the family size is a binomial with parameters $\cn$
and $p$, and the initial population size is 1.

\begin{lemma}
[Stochastic domination of cluster size by
branching process progeny size]
\label{lem-upper}
For every
$\ell \in \bbN$,
    \begin{align*}
    \Pr (|Q({\bf v}_0)|\geq \ell) \leq
    \prob_{\cn,p}(Q\geq \ell).
    \end{align*}
\end{lemma}

The second one is a slight extension of \cite[Lemma 4.3]{hl},
and establishes a lower bound. Let $\cn'= \cn - \frac{5}{2} \max
\{\ell n^{-1},C \log n\}$ and note that $\cn' \ge 2(n-1)-\frac12
\varepsilon n$ for $n$ sufficiently large. It turns out that the
cluster size can be stochastically bounded from below using a
Galton-Watson process where the family size is a binomial with
parameters $\cn'$ and $p$, and the initial population size is 1. For
$n$ sufficiently large, this process is supercritical, with mean
population size at least $1+\vep/2$, since $\cn' \ge 2(n-1)-\frac12
\varepsilon n$.

\begin{lemma}
[Stochastic domination of cluster size over
branching process progeny size]
\label{lem-lower}
There is a constant $C > 0$ such that the following holds. For every $\ell \le
\vep n^2/5$,
    \begin{align}
    \label{eqn.lower-1}
    \Pr (|Q({\bf v}_0)|\geq \ell) \geq
    \prob_{\cn',p}(Q\geq \ell)+O(n^{-6}),
    \end{align}
where $\cn'= \cn - \frac{5}{2} \max \{\ell n^{-1},C \log n\} \ge 2(n-1)-\frac12 \varepsilon n$.
\end{lemma}
Lemma~\ref{lem-lower} can be proved in exactly the same way as
\cite[Lemma 4.3]{hl}, using an
extension of \cite[Proposition 4.4]{hl} concerning the number of
elements per line in large clusters from $\eta \ll \vep$ to $\eta
\le \vep/5$ (which is exactly the same, again, since the proof of
\cite[Proposition 4.4]{hl} does not in any way rely on $\eta$
being of a smaller order than $\vep$).

In fact, Lemmas~\ref{lem-upper} and~\ref{lem-lower} are not
sufficient for our purposes, and we refine them in the following.
Let $Q_t ({\bf v_0}), G_t ({\bf v_0})$ denote the total number
of vertices and the number of unexplored vertices at time $t$ in the
exploration of the cluster of vertex ${\bf v_0}$. Also, $Q_t$ and $G_t$,
respectively, will denote the total number of nodes and the number of
unexplored nodes at time $t$ in the Galton-Watson tree when the
offspring is binomial $\Bi (\cn,p)$; and let $Q'_t$ and $G'_t$,
respectively, denote the total number of nodes and the number of
unexplored nodes at time $t$ in the Galton-Watson tree when the
offspring is binomial $\Bi (\cn',p)$.
Let ${\mathcal E}_t$
be the event that, for every $i$, no more than $m$ vertices $(i,x)$
and no more than $m$ vertices $(x,i)$ have been included in the
cluster of a vertex ${\bf v_0}$ up to time $t$ during its exploration
process. Also, let ${\mathcal E}'_t$ be the event that $Q'_t \le Q_t ({\bf v_0}) \le
Q_t$ and $G'_t \le G_t ({\bf v_0}) \le G_t$.

\begin{lemma}[Sandwiching the cluster exploration]
\label{lem.couple}
Let $\ell= \vep n^2/5$, and let $m =  5\ell n^{-1}/2$.
Then, if $n$ is large
enough, on the event
${\mathcal E}_t \cap {\mathcal E}'_t$, there exists a coupling
$\prob_{\cn, \cn',p}$ of the cluster
exploration process and the upper and lower bounding Galton-Watson
processes $\prob_{\cn,p},\prob_{\cn',p}$ in such a way that,
$\prob_{\cn, \cn',p}$-almost surely,
$Q'_{t+1} \le Q_{t+1} ({\bf v_0}) \le
Q_{t+1}$ and $G'_{t+1} \le G_{t+1} ({\bf v_0}) \le G_{t+1}$.
\end{lemma}
It is easy to prove Lemma~\ref{lem.couple} using
standard component exploration and coupling methods, in a similar
way to \cite[Lemmas 4.1 and 4.3]{hl}, and so we
omit the details.
We are now ready to prove that there is indeed no middle ground:
\medskip

\begin{proofof}{Lemma~\ref{lem.no-middle-ground}}
Lemma~\ref{lem.couple} implies that, on the event ${\mathcal E}_t
\cap {\mathcal E}'_t$, the $(t+1)^{\rm th}$ step of the exploration
process of the cluster of vertex ${\bf v_0}$ can be coupled with the
$(t+1)^{\rm th}$ step of a narrow-banded process in Lemma~\ref{lem.key}.
Since $[2(n-1)-\frac12 \vep n] p \ge 1 + \vep/2$ for $n$ large
enough, the family size of the narrow-banded branching process in
question (i.e., the component exploration process) falls into the
interval $[1+\vep/2,1+\vep]$, as required. Now observe that
$Q({\bf v_0}) \ge \ell$ if and only if  $Q_{\ell } ({\bf v_0})
\ge \ell$. We use this fact, first with $\ell = 2^{8} \vep^{-2}
\log (n^2 \vep^{-3})$, and then with $\ell = \vep^2 n/5$.
Then, the first claim follows directly from Lemma~\ref{lem.key}, 
also noting that $\prob_p ({\mathcal E}_t^c) = O(n^{-6})$ for all $t
\le \vep n^2/5$, see \cite[Proposition 4.4 and its proof]{hl}.

As for the second claim, note that for every $x, y$, the number of
components of size in between $x$ and $y$, where $0\leq x\leq y$, is given by
    \eqn{
    N_{x,y} = \sum_{{\bf v}} \frac{1}{|Q({\bf v})|} \ind_{x\leq |Q({\bf v})| < y}.
    }
Let $x = 2^8 \vep^{-2} \log (n^2\vep^3)$ and $y=\vep n^2/5$.
Then, for any vertex ${\bf
v}$,
    \begin{eqnarray}
    \E \Big (\frac{1}{|Q({\bf v})|} \ind_{x\leq |Q({\bf v})| < y} \Big ) &
    \le & \frac{2^{-8} \vep^2 }{\log (n^2\vep^3)}
    \prob_p \Big (2^8 \vep^{-2} \log (n^2\vep^3)\leq |Q({\bf
    v})|
    <\vep n^2/5\Big )\nonumber\\
    &  \le &  \frac{C2^{-8} \vep^2 }{\log (n^2\vep^3)} \Big(\vep e^{-\log
    (n^2\vep^3)}+n^{-6}\Big)\leq \frac{C}{n^2\log (n^2\vep^3)},
    \end{eqnarray}
where, for the second inequality, we have used \refeq{no-middle-ground}.
Summing over all vertices ${\bf v}$, we see that $\E [N_{x,y}]=o(1)$,
and hence $\prob_p (N_{x,y} \ge 1) =o(1)$, as required.

\end{proofof}
\smallskip

\noindent
We now complete the proof of Theorem~\ref{thm-main}:

\begin{proofof}{Theorem~\ref{thm-main}}
By Lemmas~\ref{lem.large-clusters}--\ref{lem.variance}, $Z_{\sss \geq 2^8
\vep^{-2} \log (n^2\vep^3)}$, the number of vertices in components of
size at least $2^8 \vep^{-2} \log (n^2\vep^3)$, is concentrated around $2\vep n^2$. In
other words, the number of vertices in connected components of size
at least $2^8 \vep^{-2} \log (n^2\vep^3)$ is close to $2\vep n^2$
\whp. Now, from~\cite{hl}, we know that, \whp{},
there is a giant component of size $2 \vep n^2 (1+\op (1))$. This
implies that, \whp{}, there is no other cluster of
size at least $\vep n^2/5$. Further, by \refeq{no-comp} in
Lemma~\ref{lem.no-middle-ground}, \whp{}, there are no
components of size at least $2^8 \vep^{-2} \log (n^2\vep^3)$ and less
than $\vep n^2/5$. Hence, \whp{}, the second largest component must be at
most $2^8 \vep^{-2} \log (n^2\vep^3)$, as claimed.
\end{proofof}

\subsection*{Acknowledgement}
This work was started during a visit by MJL to NYU, and continued
during a visit by MJL and JS to Georgia Tech, and during a visit by
MJL to Eurandom. The hospitality of the three institutions is
gratefully acknowledged. Also, the work of RvdH was supported in
part by Netherlands Organisation for
Scientific Research (NWO), and the work of MJL
by the Nuffield Foundation.

\newcommand\AAP{\emph{Adv. Appl. Probab.} }
\newcommand\JAP{\emph{J. Appl. Probab.} }
\newcommand\JAMS{\emph{J. \AMS} }
\newcommand\MAMS{\emph{Memoirs \AMS} }
\newcommand\PAMS{\emph{Proc. \AMS} }
\newcommand\TAMS{\emph{Trans. \AMS} }
\newcommand\AnnMS{\emph{Ann. Math. Statist.} }
\newcommand\AnnPr{\emph{Ann. Probab.} }
\newcommand\CPC{\emph{Combin. Probab. Comput.} }
\newcommand\JMAA{\emph{J. Math. Anal. Appl.} }
\newcommand\RSA{\emph{Random Struct. Alg.} }
\newcommand\SPA{\emph{Stoch. Proc. Appl.} }
\newcommand\ZW{\emph{Z. Wahrsch. Verw. Gebiete} }
\newcommand\PTRF{\emph{Probab. Theor. Relat. Fields}}
\newcommand\DMTCS{\jour{Discr. Math. Theor. Comput. Sci.} }

\newcommand\AMS{Amer. Math. Soc.}
\newcommand\Springer{Springer}
\newcommand\Wiley{Wiley}

\newcommand\vol{\textbf}
\newcommand\jour{\emph}
\newcommand\book{\emph}
\newcommand\inbook{\emph}
\def\no#1#2,{\unskip#2, no. #1,} 

\newcommand\webcite[1]{\hfil\penalty0\texttt{\def~{\~{}}#1}\hfill\hfill}

\def\nobibitem#1\par{}

\end{document}